\magnification =\magstep1
\baselineskip =13pt
\overfullrule =0pt
\centerline {\bf On the $q$-analog of homological algebra}
\vskip 1cm
\centerline {\bf M.M. Kapranov}
\vskip .5cm
\beginsection \S0 Introduction.

\vskip .5cm

 Homological algebra can be seen as the study of the equation $d^2=0$ 
(see the epigraph to [1]). It is natural to ask why $d^2$ and not, say, 
$d^3$. Following this mood, we give the natural definition.

\proclaim Definition 0.1. Let ${\cal A}$ be some Abelian category (e.g.
 the category of modules over some ring) and $N\geq 1$ be an integer.
 An $N$-complex in ${\cal A}$ is a sequence of objects and morphisms 
of ${\cal A}$
$$ C_\cdot = \{...\rightarrow C_{1}\rightarrow C_0\rightarrow C_{-1}
\rightarrow ...\}\eqno (0.1)$$
in which the composition of any $N$ consecutive morphisms equals $0$.

Thus a 1-complex is just a graded object (no differential) and a 2-complex 
is a chain complex in the usual sense.

A sourse of 2-complexes is provided by simplicial sets. Recall [1] that a
 simplicial set $X_\cdot$ is a collection of sets $X_i, i\geq 0$, 
 maps $\partial_i :X_{n}\rightarrow X_{n-1}, i= 0,1,..., n$ (called face
 maps) given for any $n$ and similar degeneracy maps 
$s_i:X_n\rightarrow X_{n+1}, i=0,1,...,n+1$. These maps are subject to 
certain commutation relations of which we note the following:
$$\partial_i\partial_j = \partial_{j-1}\partial_i  \quad {\rm for}\,\,\,
\,\,\, i<j \eqno (0.2).$$

The chain complex of a simplicial set $X_\cdot$ is the complex of 
{\bf C}- vector spaces ${\bf C}[X_\cdot]$ whose $n$- th term is
 ${\bf C}[X_n]$, the vector space freely generated by the set $X_n$. 
The standard chain differential $d:{\bf C}[X_n]\rightarrow {\bf C}[X_{n-1}]$ 
is defined by the formula $d=\sum (-1)^i \partial_i$ and satisfies $d^2=0$.
 Let $q$ be any complex number.   Define the $q$-analog of the differential
 $d$ to be the map
$$d_q:{\bf C}[X_n]\rightarrow {\bf C}[X_{n-1}],\quad 
d_q =\sum_{i=0}^n q^i\partial_i.\eqno (0.3)$$
For $q=-1$ we obtain the usual formula. 

\proclaim Proposition 0.2. Let $q$ be a $N$-th root of unity,
 $q^N=1, q\neq 1$. Then the sequence $({\bf C}[X_\cdot], d_q)$ is
 an $N$-complex.

To prove Proposition 0.2, we introduce, as it is common in the theory
 of $q$-analogs, basic numbers
$$[n]_q = (1-q)^n/(1-q) = 1 + q + ... + q^{n-1}, \eqno (0.4)$$
and basic factorials
$$[n!]_q = [1]_q. [2]_q.... [n]_q = \sum_{w\in S_n} q^{l(w)}, \eqno (0.5)$$
where in the last sum $S_n$ is the symmetric group on $n$ letters and
 $l(w)$ is the length of a permutation $w$. If $q$ is a non- trivial 
$N$-th root of unity
then $[N!]_q =0$. Therefore Proposition 0.2. is an immediate consequence 
of the following more precise statement whose proof is left to the reader.

\proclaim Lemma 0.3. For any $q\in {\bf C}$, any simplicial set $X.$, and 
simplex $x\in X_m$ and any $N\leq m$ we have 
$$d_q^N(x) = [N!]_q \sum_{i_1\geq ...\geq i_N} q^{i_1+...+i_N} \partial_{i_1}.
..\partial_{i_N}(x). \eqno (0.6)$$

The purpose of this note is to demonstrate that there is a meaningful
 homological algebra of $N$-complexes. This formalism includes $q$-
analogs of  classical  combinatorial functions taken for $q$  an $N$-
th root of unity. These values of $q$ have recently attracted much 
attention  in the study of representations of quantum groups. 
It seems that "$N$-homological algebra might be relevant to  this study. 

We show that homology objects of an $N$-complex form naturally an $(N-1)$-
complex. The role of Euler charactersitic for $N$-complexes is played by 
the value of Poincar\'e polynomial at $N$-th root of unity. In particular,
 for an exact $N$-complex this value is 0. 

In \S2 we construct the $q$-analog of the de Rham complex adding to
 commuting coordinates $x_i$ a set of $q$-commuting variables $dx_i$ 
and equip the resulting algebra with the difeferential $d$ satisfying 
the $q$-analog of Leibniz rule. For $q$ a primitive $n$-th root of 1,
 the differential satisfies the
equation $d^N=0$.  This differential calculus is not covariant under usual
 changes of coordinates, but there is a natural quadratic Hopf algebra acting
 on $x_i$ from the left and on $dx_i$ from the right, which gives "quantum 
covariance".

For $q^N=1$ and a matrix-valued 1-form $A$ ( in the $q$-deformed sense) we 
show that for  the operator $(d+A)^N$ is given by the multiplication with 
some matrix $N$-form called the  curvature of $A$.  In [4] G.Lusztig 
suggested an analogy between "quantum  geometry" at roots of unity and
 algebraic geometry in characteristic $p$. Our construction of curvature
 can be seen in this vein as the analog of the $p$-curvature of connection 
in characteristic $p$ (which measures non-commuting of covariant derivative
 with Frobenius map). In particular, for $N=3$ we obtain an expression
 resembling Chern-Simons functional in usual gauge theory.

\hfill\vfill\eject
 
\beginsection \S1. Homology of $N$-complexes.

\vskip 1cm

Let $(C., d_C), (E.,d_E)$ be two $N$-complexes in category ${\cal A}$.
 A {\it chain morphism} $f:(C., d_C) \rightarrow (E.,d_E)$ is a collection
 of morphisms 
 $f_n: C_n \rightarrow E_n$ commuting with differentials.We shall denote by 
 $N-Com ({\cal A})$ the category of all  $N$-complexes in ${\cal A}$ with 
 morphisms just defined.
\vskip .3cm
\noindent {\sl A. Homology.}

\proclaim Definition 1.1. Let $(C_\cdot, d)$ be an $N$-complex over 
an Abelian category ${\cal A}$. Its homology are the objects
$$_pH_i(C_\cdot) = {Ker \{d^p:C_i\rightarrow C_{i-p}\}\over
 Im \{d^{N-p}:C_{i+N-p}\rightarrow C_i\}}\eqno (1.1)$$
where $i\in {\bf Z}, p=1,2,...,N$. An $N$-complex is called $N$-exact 
if all its homology objects are 0.

 Clearly the homology is functorial with respect to chain morphisms of
 $N$-complexes. 

\noindent {\bf Example 1.2.} Any sequence with only $N$ consecutive 
terms
$$...0\rightarrow C_{N-1}\buildrel d_{N-1}\over\rightarrow 
...\buildrel d_2\over\rightarrow C_1\buildrel d_1\over\rightarrow
 C_0\rightarrow 0...\eqno (1.2)$$
is an $N$-complex. The condition of $N$-exactness of this $N$-complex 
means that all the $d_i$ are isomorphisms. This can be easily seen by
 induction starting from the right end of the complex. 

For $N\geq 3$ the homology of $(C_\cdot, d)$ are connected by  two
 families of natural maps. The first family is formed by
$$i_*:_pH_i(C.)\rightarrow _{p+1}H_i(C.),\quad z\,\,\,
 mod\,\,\, Im(d^{N-p})\mapsto z \,\,\,mod\,\,\, Im(d^{N-p-1})\eqno (1.3)$$
which are defined since $Ker(d^p)\i Ker (d^{p+1})$ and 
$Im(d^{N-p})\supset Im(d^{N-p-1})$. The second family consists of 
$$d_*:_pH_i(C.)\rightarrow _{p-1}H_{i-1}(C_.),\quad 
z\,\,\, mod\,\,\, Im(d^{N-p})\mapsto dz\,\,\, mod\,\,\,
 Im(d^{N-p+1}).\eqno (1.4)$$
The morphisms of the homology induced by a chain morphism 
of complexes obviously preserve these maps. We shall arrange
 $_pH_i$  pictorially in such a way that
$i_*$ acts horizontally and $d_*$- vertically. For example,
 the homology of a 5-complex look like
$$\matrix{ \downarrow&&\downarrow&&\downarrow&&&&&&&&\cr
_1H_{i+1}&\buildrel i_*\over\rightarrow &_2H_{i+1}&\buildrel i_*\over
\rightarrow &_3H_{i+1}&\buildrel i_*\over\rightarrow & _4H_{i+1}&&&&&&\cr
&&\downarrow&&\downarrow&&\downarrow &&&&&&\cr
&&_1H_i&\buildrel i_*\over\rightarrow & _2H_i&\buildrel 
i_*\over\rightarrow &_3H_i&\buildrel i_*\over\rightarrow&_4H_i &&&&\cr
&&&&\downarrow &&\downarrow &&\downarrow &&&&\cr
&&&&_1H_{i-1}&\buildrel i_*\over\rightarrow&_2H_{i-1}&\buildrel 
i_*\over\rightarrow &_3H_{i-1}&\buildrel i_*\over\rightarrow&_4H_{i-1}&&\cr
&&&&&&\downarrow&&\downarrow&&\downarrow &&\cr
&&&&&& _1H_{i-2}&\buildrel i_*\over\rightarrow&_2H_{i-2}&
\buildrel i_*\over\rightarrow&_3H_{i-2}&\buildrel i_*\over
\rightarrow&_4H_{i-2}\cr
&&&&&&&&\downarrow &&\downarrow &&\downarrow }\eqno (1.5)$$

The squares in thes diagram obviously commute. 
We form the "total object" of the homology diagram by 
setting
$${\bf H}_m(C.) = \bigoplus _{2i-p=m} {_pH_i}(C.)\eqno (1.6)$$
and define morphisms $D :{\bf H}_m(C.)\rightarrow {\bf H}_{m-1}(C.)$ 
to be equal $D = i_* +d_*$ (with the convention that $i_*$ and $d_*$
 are set to be zero when not defined. 

\proclaim Theorem 1.3. \item{(a)} Let $(C., d)$ be any $N$-complex 
in an Abelian category ${\cal A}$. Then the sequence $({\bf H}.(C.), D)$ 
is an $(N-1)$ -complex.
\item{(b)} For $N=3$, the 2-complex ${\bf H}(C.)$ is exact (in the usual
 sense).

This theorem can be reformulated by saying that for any Abelian category
 ${\cal A}$ we have a sequence of categories and functors
$$  ... \buildrel {\bf H}\over\longrightarrow 3-Com({\cal A})\buildrel {\bf H}\over\longrightarrow 2-com({\cal A})\buildrel {\bf H} \over 
\longrightarrow 1-Com({\cal A}).\eqno (1.7)$$
The part of this sequence from $3-Com ({\cal A})$ to $1-Com({\cal A})$ is 
a complex in the usual sense i.e. ${\bf H\circ H} =0$.

Before starting to prove theorem 1.3, let us consider some examples.

\noindent {\bf Example 1.4.} For a 3-complex $C.$ the (usual) exact 
sequence ${\bf H}.(C)$ has the form 
$$...\rightarrow _1H_n(C)\buildrel i_*\over\rightarrow _2H_n(C.)
\buildrel d_*\over\rightarrow _1H_{n-1}(C)\buildrel i_*\over
\rightarrow _2H_{n-1}(C)\rightarrow ... \eqno (1.8)$$

 In particular, consider a 3-complex having only two terms 
$X\buildrel f\over\rightarrow Y$. The exact sequence of homologies 
of this 3-complex has the form
$$0\rightarrow Ker(f)\rightarrow X\rightarrow Y\rightarrow Coker
 (f)\rightarrow 0\eqno (1.9)$$
More generally, any composable pair of morphisms in an Abelian category 
${\cal A}$ still gives a 3-complex $X\buildrel f\over\rightarrow
 Y\buildrel g\over\rightarrow Z$. Such a 3-complex  produces the 
six-term exact sequence
$$0\rightarrow Ker(f)\rightarrow Ker (gf)\rightarrow Ker(g)\rightarrow 
Coker(f)\rightarrow Coker(gf)\rightarrow Coker(g)\rightarrow 0.
\eqno (1.10)$$

\noindent {\bf Proof of Theorem 1.3:} Let us first prove that
 ${\bf H}.(C)$ is an $(N-1)$-complex. Since the squares in the
 homology diagram (1.5) are commutative, each matrix element of 
$D^{N-1}$ can be written as the composition of a string of morphisms 
which contains two consecutive morphisms of the form
$${_{N-1}H_i}\buildrel d_*\over\longrightarrow {_{N-2}H_i}
\buildrel i_*\over\longrightarrow {_{N-1}H_{i-1}}\eqno (1.11)$$
or
$$ _1H_i\buildrel i_*\over\longrightarrow {_2H_i}\buildrel
 d_*\over\longrightarrow {_1H_{i-1}} \eqno (1.12)$$
for some $i$. It suffices to show, therefore, that such
 compositions vanish. Indeed, to see that (1.11) gives 0,
 consider a class $[x] = x\,\,\,mod\,\,\,Im(d)\,\,\in {_{N-1}H_i}$ 
where $d^{N-1}x=0$. Then $i_*d_*[x] = dx\,\,\,mod\,\,\,Im(d)\,\,=0$.
 Similarly for (1.12). Thus ${\bf H}.(C)$ is an $(N-1)$-complex. 

Let us prove the exactness of the homology sequence (1.8) for a 3-complex. 
 First consider the term $_1H_n$. Let $w\in Ker(d)$ and 
$i_*(w \,\,mod\,\,\,Im(d^2))\,=0$. This means that $w \,\,\,
 mod\,\,\,Im(d)\,=0$ i.e. $w = dy$ for some $y\in C_{n+1}$.
 But then $y$ lies in $Ker (d^2)$ and therefore represents
 an element $[y]\in _2H_{n+1}$ so that $[w] = d_*[y]$.

Let us show the exactness  of (1.8) in the term $_2H_n$. Suppose 
that
$ d_*(z\,\,\, mod\,\, Im(d)) = 0$ i.e. $dz\,\,\,mod\,\,\,Im(d^2)\, =0$.
 This means that $dz= d^2x$ for some $x\in C_{n+1}$. Let $c = z -dx$.
 Then $dc=0$. Denoting by $[c]$ the class of $c$ in $_1H_n$, we obtain
 that $[z] = i_*[c]$ thus proving the exactness. 
\vskip .2cm
For the case of a 3-complex $C.$ Theorem 1.3.b) implies that if for
 some $p\in \{1,2\}$ all $_pH_i(C.)$ are zero then all the other
 homology groups also vanish. This fact generalizes to arbitrary
 $N$-complexes.

\proclaim Proposition 1.5. Let $C.$ be an $N$-complex such that
 for some $p\in \{1,...,N-1\}$ one has $_pH_i(C.)=0$ for all $i$.
 Then $C.$ is $N$-exact i.e. $_rH_i(C.)=0$ for all $r$ and $i$. 

\noindent {\sl Proof:} a) Consider first the case where the given 
level $p$ equals $N-1$. Let $r<N-1$. Let us show that $_rH_i=0$.
 Let $z$ be such that $d^rz=0$. We need to show that $z=d^{N-r}u$ 
for some $u$. Note that $d^{N-1}z=0$ since $r<N-1$. Therefore  
 $z=dz_1$ for some $z_1$. Further, we have $d^{N-1}z_1=d^{N-2}z=0$.
 Hence $z_1=dz_2$ for some $z_2$. Thus continuing, we find
$$z_1=dz_2, z_2=dz_3,...,z_{N-r-1}=dz_{N-r}$$
which implies $z=d^{N-r}z_{N_r}$ as required.

b) Now consider the case $r<p<N-1$. Let $d^rz=0$.  Then $d^pz=0$
 and hence $z=d^{N-p}z_1$. Since $p<N-1$, the number $N-p-1$ is
 non-negative and we can write
$$d^p(d^{N-p-1}z_1)=d^{N-1}z_1 = d^{p-1}z =0.$$
Hence $d^{N-p-1}z_1 =d^{N-p}z_2$. If $p-r>2$ then we have 
$$d^p(d^{N-p-1}z_2)= d^{p-2}z =0 \,\,\,\Longrightarrow 
\,\,\,d^{N-p-1}z_2=d^{N-p}z_3$$
for some $z_3$. Thus continuing, we find $z_3,...,z_{p-r}$
 such that $d^{N-p-1}z_i=d^{N-p}z_{i+1}$. This implies
$$z=d^{N-p}z_1=d^{N-p+1}z_2=...=d^{N-r}z_{p-r}$$
what proves the assertion.

c) Now consider the case $p<r<N-1$. Write $r=kp+l$, where 
$k,l\in {\bf Z}_+, l<p$. Let $d^rz=0$. Then $d^p(d^{r-p}z)=0$
 and so for some $z_1$ we have
$d^{r-p}z=d^{N-p}z_1$. If $k>1$ then we can continue by writing
$d^p(d^{r-2p}z-d^{N-2p}z_1)=0$ and find that  
 $d^{r-2p}z-d^{N-2p}z_1=d^{N-p}z_2$ for some $z_2$.  
  Thus continuing $k$ times, find
$$d^{r-kpz}=d^{N-kp}z_1 + d^{N-(k-1)p}z_2 + d^{N-(k-2)p}z_3 +...$$
But $r-kp =l <p$. Therefore we can write 
$d^l(z-d^{N_r}z_1 - d^{N-r+p}z_2 -...)=0$.
 By using case b), we have
$$   z-d^{N_r}z_1 - d^{N-r+p}z_2 -... = d^{N-l}w$$
which implies that $z\in Im(d^{N_r})$. The proposition is proven.               

\vskip .3cm
\noindent {\sl  B.Homotopies}
\vskip .3cm
 It is possible to define the analog of homotopies of morphisms. 
Usually it is done as a particular case of internal $Hom$-
 construction for complexes. We shall use the same approach.
 To this end we shall make the following convention.

\proclaim Convention 1.6. In the sequel we shall consider only 
categories which are {\bf C}- linear. By $\epsilon _N$ or $\epsilon$,
 when $N$ is clear from the context, we shall always denote the primitive 
root of unity $exp (2\pi i/N)$.

 For completeness sake we shall consider, along with $N$-complexes
 in a category ${\cal A}$, arbitrary sequences of the form (0.1)
 (which will be called just {\it sequences}). Morphisms of sequences 
are defined similarly to morphisms of $N$-complexes. The category of
 sequences in ${\cal A}$ will be denoted $Seq({\cal A})$.

\proclaim Definition 1.7.
\item{(a)} Let $(C.,d_C),(E.,d_E)$ be two sequences in an
 Abelian category ${\cal A}$. Their $q- Hom$- sequence is 
the sequence $\underline {Hom}(C., E.)$ of vector spaces which 
has terms
$$\underline {Hom}(C.,E.)_n = \prod_i Hom_{\cal A}(C_i, E_{i+n}).
\eqno (1.13)$$
If $(f_i:C_i\rightarrow E_{i+n})$ is an element of $\underline
 {Hom}(C.,E.)_n$ then its differential in the sequence 
$\underline {Hom}$ equals
$$d(f_i) = (g_i:C_i\rightarrow E_{i+n-1}),\,\,\,\,{\rm where}
 \,\,\,\,g_i= d_Ef_i - q^i f_{i+1}d_C.\eqno (1.14)$$
\item {(b)} Let $(V.,d_V), (W.,d_W)$ be two sequences of
 {\bf C}- vector spaces. Their $q$-tensor product is the
 sequence $V.\otimes W.$ defined by
$$(V.\otimes W.)_n =\bigoplus_{i+j=n} V_i\otimes W_j;\,\,\,\,\, 
d(v\otimes w) = d_V(v)\otimes w + q^{deg(v)}v\otimes d_W(w)\eqno (1.15)$$
where $v$ and $w$ are supposed to be homogeneous elements.

\proclaim Proposition 1.8. Let $q=\epsilon _N$ be a primitive
 $N$-th root of 1. Then the $q$-tensor product and  internal
 $q-Hom$  of two $N$-complexes  are $N$-complexes of vector spaces. 

\noindent {\sl Proof:} Let us introduce the Gaussian binomial coefficients
$$\Bigl[{N\atop k}\Bigl]_q = {(1-q^N)(1-q^{N-1})...(1-q)\over 
(1-q^k)...(1-q)(1-q^{N-k})...(1-q)}.\eqno (1.16)$$

They satisfy the recursion properties
$$\Bigl[{N+1\atop k}\Bigl]_q = \Bigl[{N\atop k}\Bigl]_q q^k + 
\Bigl[{n\atop k-1}\Bigl]_q. \eqno (1.17)$$
and are symmetric with respect to the replacement of $k$ by $N-k$. 
These properties imply the following fact.

\proclaim Lemma 1.9. Let $q\in {\bf C}$ and   $(V.,d_V), (W.,d_W)$ - 
two sequences of vector spaces, $v\in V_a, w\in W_b$. Then for any 
$n\geq 0$ we have the following equality in the $q$-tensor product
 $V.\otimes W.$
$$d^n(v\otimes w) = \sum_{k=0}^n q^{(n-k)a}\Bigl[{n\atop k}\Bigl]_{q^{-1}}
 d^ku \otimes d^{n-k}v.\eqno (1.18)$$
If $C.,E.$ are two sequences in an Abelian category ${\cal A}$ then
 for any collection $f = (f_i:C_i\rightarrow E_{i+a})$ in
 $\underline {Hom}(C.,E.)_a$ and any $n\geq 0$ we have
 the following equality in the $q- \underline{Hom}$ -sequence:
$$d^n(f) = \sum_{k=0}^n (-1)^{n-k} q^{(n-k)a} 
\Bigl[{n\atop k}\Bigl]_{q^{-1}} d_E^k \circ f 
\circ d_C^{n-k}.\eqno (1.19)$$

Proposition 1.8 follows from Lemma 1.9 when we note that for 
$q=\epsilon _N$   all the Gaussian coefficients 
$\Bigl[{N\atop k}\Bigl]_{q^{-1}}$ vanish except 
$\Bigl[{N\atop 0}\Bigl]_{q^{-1}} =
 \Bigl[{N\atop N}\Bigl]_{q^{-1}} =1$.

\proclaim Proposition 1.10. Fix a base $q\in {\bf C}$.
  For any three sequences $C.,D.,E.$ in a category ${\cal A}$ 
we have a natural morphism of sequences of vector spaces 
$$ \underline{Hom}(C.,D.)\otimes \underline{Hom}(D.,E.)
\rightarrow \underline{Hom}(C.,E.)$$
which takes an element $(f_p)\otimes (g_p)$ of 
$\underline{Hom}(C.,D.)_m\otimes \underline{Hom}(D.,E.)_n$
 to the element $(\epsilon)^{mn}g_{p+m}\circ f_p)$ of 
$\underline{Hom}(C.,E.)$.

\noindent {\sl Proof:} A straightforward checking.

Proposition 1.10. means that for the differential of the 
composition of morphisms (not necessarily chain morphisms) 
we have the $q$-Leibniz rule:
$$d(fg) = (df)g + q^{deg (f)}f (dg).\eqno (1.20)$$
\vskip .3cm

It is clear from definition 1.7 that a chain morphism 
 between $N$-complexes $C.$ and $E.$ is nothing but an
 element of degree 0 in the complex  $\underline{Hom}(C.,E.)$
 annihilated by the differential $d$. Call a morphism 
$f:C.\rightarrow E.$ {\it  null- homotopic } if it lies 
in the image of the operator $d^{N-1}$ in $\underline{Hom}(C.,E.)$.
 In other words, $f$ is null - homotopic  if there exists a
 collection of morphisms $s_i: C_i\rightarrow E_{i+N-1}$
 such that 
$$f_i = s_{i-N+1}d^{N-1} +\epsilon d s_{i-N+2} d^{N-2} + ... + 
\epsilon ^{N-1} d^{N-1}s_i.$$
This follows from Lemma 1.9 after  suitable renormalization of  $s_i$. 

\proclaim Proposition 1.11. A null- homotopic morphism 
 $f:C.\rightarrow E.$ of $N$-complexes induces a zero map
 on all the homology objects $_pH_i$.

\noindent {\sl Proof:} Let $z\in C_i$ and $d^pz=0$ for some $p<N$.
 Then 
$$f(z) = \epsilon ^{N-p} d^{N-p} s_{i-N+1} d^{p-1}z +\epsilon ^{N-p+1}
 d^{N-p+1}s_{i-N+2} d^{p-2}z +... + \epsilon ^{N-1}d^{N-1}s_iz.$$
This expression obviously lies in $Im(d^{N-p})$ so the class of $f(z)$
 in $_pH_i$ equal 0, QED.

\proclaim Proposition 1.12. Null- homotopic morphisms form an ideal in the
 category $N-Com ({\cal A})$.

\noindent {\sl Proof:} Let $f,g$ be two morphisms of $N$-complexes 
such that the composition $fg$ is defined. Suppose that $f$ is
 null- homotopic, i.e. $f= d^{N-1}s$. Since $dg=0$, we have, by
 iterating Leibniz rule (1.9) that $d^{N-1}(sg)= (d^{N-1}s)g =fg$
 so $fg$ is null- homotopic. Similarly we prove that if $g$ is 
null- homotopic that so is $fg$.
\vskip .3cm

The quotient category of  $N-Com({\cal A})$ modulo the ideal
 of null- homotopic morphisms will be denoted $N- Hot({\cal A})$
 and called the {\it homotopy category of $N$-complexes}. 

\vskip .3cm
\noindent {\sl C. $q$-Euler charactersitic.}
\vskip .3cm
We shall consider, for simplicity, finite $N$-complexes of
 finite- dimensional complex vector spaces. If $C.$ is such
 an $N$-complex, we define its Poincar\'e polynomial
 $P_C(q) = \sum dim (C_i) q^i$. The usual Euler characteristic 
 for a 2-complex is the value of the Poincar\'e polynomial at (-1). 
Its main property is that it vanishes for an exact 2-complex. 
Let us give a generalization of this fact for $N$-complexes.

\proclaim Proposition 1.13. Let $C.$ be an $N$-exact $N$-complex 
and $\epsilon$ be a primitive $N$-th root of 1. Then $P_C(\epsilon) =0$. 

\noindent {\sl Proof:} Denote individual components of the
 differential by $d_i:C_i\rightarrow C_{i-1}$. Let us denote
 the dimension of a vector space $V$ shortly $[V]$. Then we have
$$[C_i] = [Ker(d_i)] + [Im(d_i)],\eqno (1.21)$$
$$[C_i] = [Ker (d_{i+N-1}\circ ...\circ d_{i+1}\circ d_i)] + 
[Im (d_{i+N-1}\circ ...\circ d_{i+1}\circ d_i)].\eqno (1.22)$$
By $N$-exactness we have 
$$Ker(d_i) = Im(d_{i-1}\circ ...\circ d_{i-N+1}),\,\,\, 
Im(d_i) = Ker (d_{i+N}\circ ...\circ d_{i+1}). \eqno (1.23)$$
Therefore by (1.21) and (1.23) we have 
$$P_c(\epsilon) = \sum \epsilon^i([Ker (d_i)]+[Im(d_i)])=$$
 $$= \sum \epsilon^i([Im (d_{i+N-1}\circ ...\circ d_{i+1}\circ d_i)]+
[Ker (d_{i+N-1}\circ ...\circ d_{i+1}\circ d_i)])$$
Now formula (1.22) implies that $P_C(\epsilon)=\epsilon 
P_C(\epsilon)$ hence $P_C(\epsilon)=0$. 

\hfill\vfill\eject

\beginsection \S2. Differential forms and de Rham $N$-complexes.

\vskip 1cm
The theory of differential forms which we are about to develop,
 will be explicitly coordinate-dependent. Fix a complex number $q$.

Let $(x_1,...,x_n)$ be coordinates in a real affine space of dimension 
$n$. Introduce formal differentials $dx_1,...,dx_n$ which we will make 
subjects to  the relations
$$f(x) dx_i = dx_i f(x)\eqno (2.1)$$
for any differentiable function $f(x)$ and
$$dx_idx_i =0;\,\,\,  dx_i dx_j = q dx_j dx_i \,\,\,\, {\rm 
for}\,\,\,\, i>j \eqno (2.2)$$
The algebra generated by functions $f(x)$ and the symbols $dx_i$ 
will be called the ($q$-analog of) algebra of differential forms 
on ${\bf R}^n$ and denoted $\Omega^\cdot_{q,{\bf R}^n}$. It has
 the obvious grading where functions have degree 0 and $dx_i$
 have degree 1. A homogeneous element of degree $k$ (a $k$- form)
 can be written uniquely in the form
$$\omega =\sum_{1\leq i_1<...<i_k\leq n}f_{i_1,...,i_k}(x) 
dx_{i_1}...dx_{i_k} \eqno (2.3)$$
and thus can be identified with the collection of $f_{i_1,...,i_k}(x)$. 
Define the exterior differential of the form (2.3) to 
$$d\omega =\sum_{1\leq j_1<...<j_{k+1}\leq n}g_{j_1,...,j_{k+1}}(x)
 dx_{j_1}...dx_{j_{k+1}},\eqno (2.4)$$
where
$$g_{i_1,...,i_{k+1}}(x)= \sum _{p=1}^{k+1}q^{p-1}
 {\partial f_{j_1,...,\hat j_p,...,j_{k+1}}\over \partial x_p}
\eqno (2.5)$$

\proclaim Proposition 2.1. 
\item{(a)} The differential $d$ satisfies the $q$-analog of
 the Leibniz rule: for $u\in\Omega^i, v\in\Omega^j$ we have
$$d(uv) = d(u)v + q^i u d(v).\eqno (2.6)$$
\item{(b)} Suppose that $q$ is a primitive $N$-th root of unity.
 Then $d^N =0$.

\noindent {\sl Proof:} For 0-forms the differential is given by 
the usual formula: $df = \sum (\partial f/\partial x_i) dx_i$.
 The formula (2.5) expresses exactly the $q$-Leibniz rule for
 $d(f(x).dx_{i_1}...dx_{i_k})$. The $q$-Leibniz rule for a 
general product follows from this. Part (b) of Proposition 2.1 
is a consequence of the following lemma which we leave to the reader.

\proclaim Lemma 2.2. Let $A=\bigoplus A_i$ be a graded ${\bf C}$- 
algebra and $d:A\rightarrow A$ be a linear map of degree +1 satisfying 
the $q$-Leibniz rule (2.6) for some $q\in {\bf C}$. Then for each 
$N\geq 0,\, u\in A_i,v\in A_j$ we have 
 $$d^N(uv) = \sum _{p=0}^N q^{ip}\Bigl[{N\atop p}\Bigl]_q d^p(u) 
d^{N-p}(v)\eqno (2.7)$$

\noindent {\bf Remark 2.3.} Our $q$-deformed de Rham differential 
can be seen as a particular case of the construction of \S0, which
 produces an $N$-complex out of every simplicial vector space. 
Indeed, our $\Omega^k_q$ as vector spaces (and modules over functions)
 do not depend on $q$, though their exterior multiplication does.
 Define vector spaces $V_{-1} = \Omega^n, V_0 = \Omega^{n-1}$ etc. 
For any $p$ introduce linear operators $\partial_\nu :V_p\rightarrow 
V_{p-1}$ for $0\leq \nu\leq p$. The vector space $V_p = \Omega^{n-p-1}$ 
is a direct sum of the spaces of forms of the type 
$f(x) dx_{j_1}...dx_{j_{n-p-1}}$ for all sequences
 $1\leq j_1<...<j_{n-p-1}\leq n$. We define $\partial _i$ on
 these spaces separately. More precisely, let $i_0<...<i_p$
 be all the elements of $\{1,...,n\}$ not lying in 
$\{j_1,...,j_{n-p-1}\}$, taken in the increasing order.
 We set 
$$\partial_\nu (f(x) dx_{j_1}... dx_{j_{n-p-1}}) = ({\partial f\over 
\partial x_{i_\nu}}) dx_{j_1}... dx_{i_\nu}... dx_{j_{n-p-1}}$$
where $dx_{i_\nu}$ is inserted to such place that the indices of 
differentials increase. It is immediate to verify that operators 
$\partial_\nu$ thus defined satisfy the equations (0.2) i.e. that 
$V.$ forms an (augmented, because of $V_{-1}$) simplicial vector 
space without degenerations. Now the standard de Rham differential
 equals $\sum (-1)^\nu \partial_\nu$, whereas our $q$-deformed 
differential equals $\sum q^\nu \partial_\nu$.

\proclaim Proposition 2.4. Let $N\leq n+1$. Then the cohomology 
spaces $_pH^i(\Omega^\cdot)$ of the de Rham $N$-complex of 
polynomial $N$-forms  equal 0 for $i+p+1 \leq N$.

\noindent {\sl Proof:} Denote the variables $dx_i$ by $\xi_i$.
 Introduce the $q$-derivations $\partial/\partial \xi_i$ on
 $\Omega^\cdot$ such that 
$$\partial/\partial \xi_i (\xi_j) = \delta_{ij},\,\, \partial/
\partial \xi_i (x_j)=0 $$
Similarly introduce usual derivations $\partial/\partial x_i$.
Then the exterior differential $d$ equals $\sum \xi_i\partial/
\partial x_i$.
Consider the operator $S= (\sum x_i\partial/\partial \xi_i)^{N_1}$
 as a homotopy in the de Rham complex. We introduce in
 $\Omega^\cdot$ the double grading $\bigoplus \Omega^i(j)$
 where $\Omega^i(j)$ is the subspace of forms homogeneous
 of degree $i$ with respect to $\xi$'s and of degree $j$ 
with respect to $x$'s.

The map $d^iSd^{N_1-i}$ is calculated explicitly to be an 
isomorphism on each $\Omega ^r(s)$ with $r+s \leq N-i$.
 The assertion follows from this.

\vskip .3cm
 Let us now discuss the question of covariance of the constructed
 exterior differential calculus. Since the construction makes an 
explicit appeal to coordinates, our algebra of forms and the exterior
 differential are not invariant under the group $GL(n)$ of linear 
changes of coordinates. However, it is possible to introduce a certain
 Hopf algebra which preserves all the picture. This algebra is a 
version of quantized $GL(n)$ which we shall now describe.

Let $ a_{ij}, \,\,i,j=1,...,n$ be independent non- commutative
 generators which we arrange into a matrix $A=||a_{ij}||$. We multiply
 formally $A$ from the left to the  row vector of variables
 $x=(x_1,...,x_n)$, obtaining $Ax$, where $(Ax)_i = \sum a_{ij}x_j$.
 Similarly we multiply $A$ from the left to the column vector 
 $dx = (dx_1,...,dx_n)^t$, obtaining $dx.A$, where 
$(dx.A)_i = \sum A_{ji}dx_j$.
Let us now require that $(Ax)_i$ and $(dx.A)_i$ commute 
in the same way as $x_i$ and  $dx_i$ do i.e.:
$$(Ax)_i(Ax)_j = (Ax)_j(Ax)_i \,\,\,\forall \,\,i,j\,\,{\rm and}
\,\,\, (dx.A)_i(dx.A)_j = $$
$$= q (dx.A)_j(dx.A)_i \,\,\,{\rm for}\,\,\,\, i> j,\,\,(dx.A)_i^2=0.
\eqno (2.8)$$

Here we understand that $x_i$ or $dx_i$ satisfy the imposed relations.
 By taking coefficients of relations (2.8) by $x_ix_j$ or $dx_jdx_j,
 i\leq j$, we obtain explicit relations:
$$a_{ij}a_{ik}=a_{ik}a_{ij};\,\, a_{ij}a_{kj}=q a_{kj}a_{ij}\eqno (2.9)$$
$$a_{ii}a_{jj}-a_{jj}a_{ii} = a_{ji}a_{ij}- a_{ij}a_{ji}= 
q^{-1}a_{ij}a_{ji}-qa_{ji}a_{ij}\,\,\, {\rm for}\,\,\,  i<j\eqno (2.10)$$

Denote by $R(n)$ the associative algebra generated by
 $a_{ij}, \,i,j\leq n$ subject to the relations (2.9) and (2.10).

\proclaim Proposition 2.4. \item{(a)} The  formula $\Delta(a_{ij}) = 
\sum_k a_{ik}\otimes a_{kj}$ defines a Hopf algebra structure 
$\Delta :R(n)\otimes R(n)\rightarrow R(n)$ on $R(n)$.
\item{(b)} The formulas $\alpha (x_j) = \sum_i a_{ij}\otimes x_i,
 \alpha (dx_i) = \sum _j a_{ij}\otimes dx_j$ defines a lefy coaction 
$\alpha :\Omega^\cdot_q\rightarrow R(n)\otimes \Omega^\cdot_q$ of $R(n)$
 on the $q$-de Rham algebra.
\item{(c)} The coaction $\alpha$ commutes with the $q$-exterior
 differential $d$ i.e., $\alpha\circ d = (1\otimes d)\circ\alpha$.

\noindent {\sl Proof:} Our construction of the algebra $R(n)$ 
is a slight modification of the construction of internal $Hom$ 
of quadratic algebras introduced  by Y.I.Manin in [2].
 Let $A$ be a homogeneous quadratic algebra  with generators 
and relations $(V, I\i V\otimes V)$. For such $A$ Manin 
introduces a Hopf algebra $\underline{End}(A)$ of "endomorphisms"
 of these relations so that $\underline{End}(A)$ coacts on $A$
 from the right. Similarly one can consider another group of
 relations on the dual space $(V^*,J\i V^*\otimes V^*)$ and 
require a left coaction on the quadratic algebra $D$ defined
 by latter relations. One can require both coactions, thus 
obtaining a Hopf algebra $\underline{End}(D,A)$. It is
 straightforward to see, using methods of [2] that 
$\underline{End}(D,A)$ is always a Hopf algebra.
 Our construction is a particular case of $\underline{End}(D,A)$,
 where $A$ is the symmetric algebra and $D$ is the $q$-exterior
 algebra. This establishes (a) and (b). 

\noindent {\bf Remark.} The standard $q$-analog of group
 $GL(n)$ is obtained as $\underline{End}(D,A)$ where $D$ 
is the $q$-exterior algebra, as in our case, whereas $A$ 
is a $q$-symmetric algebra $x_ix_j=qx_jx_i$, see [2].

\hfill\vfill\eject

\beginsection \S3. Connections and curvature.

\vskip 1cm
Let $q\in {\bf C}$ be a fixed complex number.
Consider a trivial vector bundle $E$ of rank $r$ over some
 domain $U\i {\bf R}^n$. Let $x_1,...,x_n$ be the standard
 coordinates in ${\bf R}^n$  and $\Omega^\cdot = \Omega^\cdot_q$- 
the graded algebra of differential forms on $U$, introduced in \S2. 
The differential in this algebra will be denoted by $d$. We denote
 by $\Omega^\cdot(E)$ the graded $\Omega^\cdot$- module formed by
 forms with values in ${\bf C}^r$. By $\Omega^\cdot (End(E))$ we
 denote the algebra $\Omega^\cdot \otimes End (E)$. The differential 
$d$ extends to this algebra and satisfies the $q$-Leibniz rule 2.4.

A {\it  $q$-connection} in $E$ is, be definition, the operator in 
$\Omega^\cdot (E)$ of the form $\nabla = d+A$, where 
$A = \sum A_i(x) dx_i$ is a matrix- valued 1-form. Such an 
operator also satisfies the $q$-Leibniz rule in the form
$$\nabla (\omega . \Sigma) = (d\omega).\Sigma + 
q^{deg (\omega)}\nabla (\Sigma),\eqno (3.1)$$
where $\omega$ is a scalar form of some degree $d= deg(\omega)$,
 and $\Sigma$ is a vector- valued form.

\proclaim Theorem 3.1. Let $q = \epsilon _n $ be a primitive
 $N$-th root of unity. Then, for any $q$-connection $\nabla$ 
the operator $\nabla^N$ is given by multiplication with some matrix- 
valued $N$-form $F_\nabla$.

We call $F_\nabla$ the {\it N- curvature} of the connection $\nabla$.

 Theorem 3.1 is a consequence of  the following lemma which 
is similar to Lemma 1.5.

\proclaim Lemma 3.2. Let $q$ be any complex number, $\nabla$ be
 any $q$-connection, $f$- any scalar function, $s$- any vector
 function. Then, for any $n\geq 0$ we have
$$\nabla^n (f.s) = \sum_{k=0}^n \Bigl[{n\atop k}\Bigl]_q (d^kf).
 (\nabla ^{n-k} s). \eqno (3.2.)$$

 Let us  prove Theorem 3.1. If $q$ is a primitive $N$ -root of 1 
 then , by (3.2), all the Gaussian coefficients
 $\Bigl[{N\atop k}\Bigl]_q$ vanish except $\Bigl[{N\atop 0}\Bigl]_q 
= \Bigl[{N\atop N}\Bigl]_q =1$. Substituting $n=N$ in formula (3.4),
 we find that 
$$\nabla^N (f.s) = (d^Nf).s + f.(\nabla ^N s) = f.(\nabla ^N s)$$
 since $d^N=0$. This means that the differential operator
 $\nabla ^N$ has in fact order 0 and thus has the claimed form.

Denote by ${\cal G}$ the group of invertible ($r\times r$)
 -matrix- valued functions on $U$ with differentiable entries.
 This group acts on connections by conjugation:
$$\nabla \mapsto g^{-1}\nabla g = g^{-1}(d+A)g.\eqno (3.5)$$
Therefore the transformation of $A$ under gauge transformation 
is given by the usual formula:
$$A\mapsto g^{-1}dg + g^{-1}A g.\eqno (3.6)$$

Clearly the curvature is covariant inder gauge transformations 
$$F_{g^{-1}\nabla g} = g^{-1}F_\nabla g.\eqno (3.7)$$

Therefore the set of connections with vanishing $N$-curvature is 
invariant under the gauge group. Any such connection $\nabla$
 produces a vector- valued de Rham $N$-complex. However, it may 
be impossible to transform, even locally, a connection with vanishing 
curvature to a trivial connection ($A=0$) by a gauge transformation. 

It is possible to prove the analog of Bianchi identity for curvatures 
and to develop a version of theory of Chern forms and Chern
 classes for $q$-connections. We shall assume for the rest
 of the section that $q=\epsilon = exp (2\pi i/N)$.

\proclaim Proposition 3.3. We have $\nabla (F_\nabla)=0$ for any $N$-connection $\nabla$.

\noindent {\sl Proof:} By the $q$-Leibniz rule we have 
$$\nabla (F_\nabla) = \nabla\circ F_\nabla -F_\nabla 
\circ \nabla = \nabla^{N+1} -\nabla^{N+1} =0.$$

\proclaim Corollary 3.4. The  scalar differential form 
$tr(F_\nabla^p)$ is closed for each $p$. 

\vskip .3cm

Consider the particular case $N=3$. Then the curvature 
of a connection $\nabla = d+A$ is the 3-form  
$$F_A = d^2(A) + d(A).A + \epsilon A.d(A) + A.A.A .\eqno (3.8)$$

The last three terms in the RHS of (3.8) have a striking 
recemblance with the famous Chern -Simons Lagrangian [3]  
(for connections in the usual sense over 3- manifolds). 
Moreover, the term $d^2A$ will note contribute to overall
 integral. In other words, we have the following fact. 

\proclaim Proposition 3.5. Let $\nabla = d+A$ be a 3-connection 
on ${\bf R}^3$ such that the functions $A_i(x)$ have compact support (or are rapidly decreasing). Then we have
$$\int_{{\bf R}^3} tr F_A = \int_{{\bf R}^3} tr ( d(A).A + 
\epsilon A.d(A) + A.A.A) \eqno (3.9)$$
where $\epsilon = exp (2\pi i/3)$.

\noindent {\sl Proof:} It suffices to show that $d^2A$ has
 zero integral. But in explicit form we have
$$d^2 A = \Bigl({\partial ^2A_3\over \partial x_1\partial x_2} 
+ \epsilon {\partial ^2A_2\over \partial x_1\partial x_3} +
 \epsilon ^2 {\partial ^2A_1\over \partial x_2\partial x_3}
\Bigl) dx_1dx_2dx_3 $$
which is a divergence.

\hfill\vfill\eject

\centerline {\bf References:}

1. S.I.Gelfand, Y.I. Manin. Methods of homological algebra. 
Part I: Introduction to cohomology theory and derived categories, 
Nauka, Moscow, 1990 (in Russian, English translation to appear in Springer).

2. Y.I. Manin. Some remarks on Koszul algebras and quantum groups, 
{\it Ann. Inst. Fourier}, {\bf 1987}, v.37, No.4, p.191-205.

3. M.F.Atyah. The geometry and physics of knots, Cambridge University 
Press, 1990. 

4. G. Lusztig, Finite dimensional Hopf algebras arizing from quantized 
universal enveloping algebras, {\it Jour. Amer. Math. Soc.}, {\bf 3} 
(1990), 257-296.

\vskip 1cm
{\sl Department of Mathematics, Northwestern University, Evanston Il 60208 USA,
email: kapranov@math.nwu.edu}

\bye